# Antinomies of Mathematical Reason:
# The Inconsistency of PM Arithmetic and Related Systems

(S. Fennell, Cambridge)

**Abstract:** We give a proof of the inconsistency of PM arithmetic, classical set theory and related systems, incidentally exposing an error in Gödel's own proof of Gödel's Theorems. The inconsistency proof, that formulae of the form R and ~R occur as theorems in the PM-isomorphic system P, proceeds from a reflexive substitution instance of the axiom $p \vee p \to p$, the first axiom of the propositional calculus (axiom II.1 of P). Gödel's formalism is used throughout.

It has long been thought that the paradox of the Cretan Epimenides uttering *Κρῆτες ἀεὶ ψεῦσται* ('Cretans are always liars'), as recorded in one of the more exiguous books of the New Testament,[1] and more earnestly expounded as a paralogism by Eubulides,[2] cannot be constructively formulated in conventional logical and mathematical systems. We formulate in the following pages, however, a structure which possesses some of the critical properties of that paralogism.

We note initially that deductive systems as a class (including mathematical systems) by their nature presuppose the relations and axioms embodied in the propositional calculus. In particular we note the indispensability of the propositional axiom (known as the 'principle of tautology') listed as axiom II.1 in Gödel's 1931 paper, $p \vee p \to p$,[3] which is equivalent to $(\sim p \to p) \to p$. If this principle, itself a tautology, were neither stipulated nor derivable in the system, little sense could be made of the system's basic connectives, and of deductive implication in particular: the consequences of excluding this axiom (on 'intuitive' grounds) would be (i) denial that $p$ is implied by $p$, and (ii) denial that $p$ is implied by $p \wedge p$ (both denials in obvious contradiction of the meanings of those signs).[4]

We shall see that cases of both R and ~R occur as theorems in the PM-isomorphic system P of basic arithmetic (where R is some well-formed formula of the system), as they will in any impredicative deductive system of this type (and in any wider set of methods or computer program[5] embracing any of those systems, as well as all theories in which the structure of the present proof can be reproduced): hence these systems are inconsistent. We shall use the same system as Gödel did in 1931; in the explanations however, for clarity's sake, we shall use 'Gn' as an abbreviation for 'Gödel-number' (where Gödel used italics to denote the Gn of any described object.) The axiom sets of the system are well known, and we need here only remind the reader of a few items which feature in the immediate formulation of Gödel's crucial self-referential formula.[6] A relationmark is a formula of P with two or more free variables. A classmark is a formula of P with a single free variable. A numbermark is the term in P corresponding to a natural number. The expressions 'Bew (*a*)', 'Neg *a*', and '*a* Imp *b*' mean: 'the P-formula with Gn *a* is provable in P', 'the Gn of the negation of the P-formula with Gn *a*', and 'the Gn of the P-formula denoting the implication of the P-formula with Gn *b* by the P-formula with Gn *a*' respectively.[7] ω-consistency is defined by Gödel as asserting that, regardless of whether the generalisation of a classmark with Gn *a* is a theorem of system P (or of some class κ of its formulae), it cannot be so that both the negation of that generalisation be a theorem, and yet each successive natural-numbermark-substitution instance of the classmark with Gn *a* also be a theorem.[8] This purports to define a stricter standard of consistency than normal consistency.

The framework for deriving contradictories in P is created essentially by expressing a variant of axiom II.1 as a formula referring to recursive self-referential formulae.

We note the recursiveness of the substitution function *Sb* as defined by Gödel, whereby the relationmark $Sb(x^{23}_{Z(y)})$ signifies the Gn of the result of substituting the numbermark of a natural-number value *y* into the single free variable (here defined as the variable with Gn 23) of the classmark with Gn *x*. We can also make use, as Gödel did, of the 'self-substitution' case $Sb(x^{23}_{Z(x)})$ where *y* = *x*: where the numbermark of a classmark's own Gn is substituted into its own single free variable. (The primary purpose of this structure is to overcome the problems of arithmetical impossibility involved in more direct strategies of self-reference.) We define a classmark with Gn *s* as follows:

$$s = Sb(x^{23}_{Z(x)}) \qquad (1)$$

Using two instances of this 'self-substitution' function, we define a classmark with Gn *t* as follows:

$$t = s \text{ Imp Neg } s \qquad (2)$$

noting that Imp is recursively defined as in Gödel's Recursive Definition 32.[9]

We now define a formula with Gn *u* as follows:

$$u = Sb\,(t^{23}_{Z(t)}) \qquad (3)$$

noting that the single free variable with Gn 23 in the classmark with Gn *t* is *x* itself. The semantic content of formula with Gn *u* is, to express it demotically: "This quoted statement implies the negation of this quoted statement." Intuitively, we might not expect this formulation to be a theorem of P, and indeed, since this formula is itself of the form ($p \rightarrow \sim p$), it follows via the P-axiom II.1, in the form ($p \rightarrow \sim p) \rightarrow \sim p$, that:

$$\text{Bew } (u \text{ Imp Neg } u) \qquad \text{[by axiom II.1]} \qquad (4)$$

for the formula *u* is also the specification of its own consequence via axiom. But since this provable implication is itself the formula referred to by its own protasis (the lefthand side of the implication formula) – i.e. it holds, in terms of Gödel's Definition 35, that $A_1$ - $Ax\,(u)$, – it follows that:

$$\text{Bew } (u) \qquad \text{[by 2, 3 \& 4]} \qquad (5)$$

Proceeding from line 4, however, the apodosis will also be a theorem of P, via the same axiom:

$$\text{Bew (Neg } u) \qquad \text{[by 4 \& axiom II.1]} \qquad (6)$$

$$\text{Bew } (u \text{ Con Neg } u) \qquad \text{[by 5 \& 6]} \qquad (7)$$

☐

This establishes that (*contra* Gödel p. 189) there is proof in P of a pair of mutually contradictory formulae of P, regardless of whether we suppose the system to be complete or not: the system P itself is inconsistent, and all speculations arising from the notion of its incompleteness must be discarded.

One immediate consequence of P's inconsistency is that the distinction between the consistent and ω-consistent formulae of this system, substantiated by formulae derived from Gödel's undecidability proof, vanishes;[10] there will exist in this system no consistent (and *a potiori* no ω-consistent) class κ of formulae in which the formulability of Gödel's undecidable formulae could be demonstrated; his 1931 proof of P's undecidability thus lapses from Theorem VI onwards.

It follows from hitherto proven elements of Church's Thesis that the same contradiction as demonstrated above arises in Turing-, lambda- and abacus-computable as well as combinatorically definable functions, Shepherd/Sturgis register machines, Markov algorithms and the geometries and other systems isomorphic to or reducible to them, or even (as extensions) simply including them. Hence the normal conditions defining these systems define in each case a contradictory complex. As the consequences of Gödel's proof would also have applied to Zermelo-Fraenkel-Neumann set theory (ZF),[11] there is no difficulty in constructing an analogous proof for ZF's inconsistency. The same applies to a vast class of theories whose supposed undecidabilites arise from self-referential formulae (or derivatives of them) which can be dealt with in the same way as above; in particular this affects all arithmetisable systems in which Cantorian diagonalisation can be defined and carried out, and thus affects the logical and mathematical results that have been obtained in that manner. Although the discussion of them goes beyond the remit of this paper, the consequences of mathematical inconsistency for systems of theoretical physics may be traced in terms of the system relations discussed in Deutsch (1985), Wolfram (1985), A.C. Manoharan (2000) and others.[12]

It is a consequence of a system's inconsistency that every proposition formulable within it (including each proposition's negation) also becomes provable; the condition for the consistency of a system is namely that the system contain at least one proposition which is not a theorem. As the computable analogues of the theorems of any such system, a Turing machine would simply give (in any formally specified order) the analogues of every 'correct' proposition and its negation: thus there are in these systems no genuinely undecidable propositions or numbering problems arising from them.

For this reason no formulation of the inconsistency principle is offered here in any such system itself: trivially, the negation of the respective inconsistency result would also be a theorem of each of these systems. The inconsistency result is in that particular sense suprasystemic. This further – and vastly more serious – disappointment for Hilbert's formalist programme changes considerably the means by which we must pursue our intuitive conceptions of verifiable mathematics and other constructive systems, although the basic impulses whose elaboration is our ongoing task remain: 'the impression that mathematics deals with structures of abstract objects that are independent of us, and the conviction, that principles for some structures are immediately evident, because we can grasp the build-up or construction of their elements'.[13]